\documentclass{article}
\usepackage{graphicx}

\newtheorem{theorem}{Theorem}

\newtheorem{prop}[theorem]{Proposition}
\newtheorem{cor}[theorem]{Corollary}
\newtheorem{definition}[theorem]{Definition}

\title{Frameworks with crystallographic symmetry}
\author{Ciprian S. Borcea and Ileana Streinu}
\date{}

\begin{document}
\maketitle

\begin{abstract}
Periodic frameworks with crystallographic symmetry are investigated from the perspective of
a general deformation theory of periodic bar-and-joint structures in $R^d$. It is shown that
natural parametrizations provide affine section descriptions for families of frameworks with a specified graph and symmetry. 
A simple geometric setting for diaplacive phase transitions is obtained.
Upper bounds are derived for the number of realizations of minimally rigid periodic graphs.
\end{abstract}



\section*{Introduction}

The notions of periodic graph and periodic framework emerged as abstractions of crystal structures. Inquiries about lattice sphere packings sprang from the same source. 
 
\medskip 
We show that a classical perspective used in the theory of positive definite quadratic forms and lattice sphere packins leads to natural parametrizations for placements of periodic graphs.
In this setting, all placements with a specified crystallographic symmetry correspond with a
certain affine section. As a result, the deformation theory for periodic frameworks presented
in our papers \cite{BS2, BS3} is adapted in a natural way so as to encompass all cases of
higher crystallographic symmetry. 

\section{Symmetries of a periodic framework}

We adopt here definitions introduced in \cite{BS2}. Let $(G,\Gamma)$ be a $d$-periodic graph. The infinite graph $G=(V,E)$ is assumed connected and when given a periodic placement $(p,\pi)$ in $R^d$, the corresponding periodic framework is denoted $(G,\Gamma,p,\pi)$. Recall that
$\Gamma \subset Aut(G)$ is a free Abelian group of rank $d$ and $\pi$ is a faithful 
representation of $\Gamma$ by a lattice of translations of rank $d$. Moreover, 

\begin{equation}\label{placement} p:V\rightarrow R^d \ \ \ \mbox{and} \ \ \ \pi:\Gamma\rightarrow {\cal T}(R^d) 
\end{equation}

\noindent
are related by:

\begin{equation}\label{relation}
 p\circ \gamma = \pi(\gamma)\circ p \ \ \ \mbox{for all} \ \gamma \in \Gamma 
\end{equation}

\medskip \noindent
Relation (\ref{relation}) shows that $\pi$ may be inferred from $p$, but most considerations
about framework deformations and symmetries benefit from observing both functions. The
quotient multigraph $G/\Gamma$ is assumed to be finite and we put 

\begin{equation}\label{nm}
n=|V/\Gamma| \ \ \ \mbox{and}\ \ \ m=|E/\Gamma|
\end{equation}

\medskip \noindent
Periodic frameworks are abstract, idealized versions of crystalline materials and,
like them, may posess other symmetries, besides those expressing periodicity under $\Gamma$.
Thus, there might be a larger group of automorphisms $\Gamma\subset \Sigma \subset Aut(G)$
and an extension of $\pi$ to a faithful representation of $\Sigma$ by a crystallographic group
$\pi(\Sigma)\subset E(d)$, such that relation (\ref{relation}) would hold for all $\sigma \in \Sigma$. 

\medskip
Considering that $Aut(G,\Gamma)$ is the normalizer of $\Gamma$ in $Aut(G)$,
a natural assumption for investigating this setup will be that $\Gamma$ is {\em normal} in $\Sigma$, that is $\Gamma \subset \Sigma \subset Aut(G,\Gamma)$. If all
translational symmetries of the framework $(G,\Gamma,p,\pi)$ are in $\pi(\Gamma)$, this
is necessarily the case, since the subgroup of translations in a crystallographic group is
normal. In general, the normality assumption would hold after replacing
the initial periodicity group $\Gamma$ by an appropriate subgroup of finite index $\tilde{\Gamma}\subset \Gamma$. Alternatively, instead of relaxing, one may refine the periodicity group by adopting all translational symmetries of the given framework. 

\medskip
For these reasons, we proceed below with the study of
framework symmetries which correspond to graph automorphisms in the normalizer 
 $N(\Gamma)$  of $\Gamma$ in $Aut(G)$. Note that the quotient group $N(\Gamma)/\Gamma$
acts naturally on the quotient graph $G/\Gamma$. It follows that $N(\Gamma)/\Gamma$ is
finite, since $G/\Gamma$ is finite and $G$ connected.

\begin{definition}
We say that $\sigma\in N(\Gamma)=Aut(G,\Gamma)$ is a symmetry of the $d$-periodic framework $(G,\Gamma,p,\pi)$ when the result
of acting by $\sigma$ on the framework is the same as the result of acting by an isometry $s\in E(d)$, that is:

\begin{equation}\label{sym}
s\circ p=p\circ \sigma 
\end{equation}

\noindent
In other words, we have a commutative diagram

\begin{equation}\label{diagram}
 \begin{array}{lcll}
                         V & \stackrel{p} \longrightarrow & R^d \\
                        \downarrow \sigma & \ & \downarrow  s \\
                         V & \stackrel{p} \longrightarrow & R^d
     \end{array}
\end{equation}
\end{definition}

\medskip \noindent
As remarked above, it is convenient to keep $p$ and $\pi$ on an equal footing. 
Then (\ref{sym}) becomes the equivalent, but more revealing condition:

\begin{equation}\label{symPair}
s\circ p=p\circ \sigma \ \ \ \mbox{and} \ \ \ C_s\circ \pi=\pi\circ C_{\sigma}
\end{equation}

\medskip \noindent
where $C_s$ denotes the restriction to the group of translations ${\cal T}(R^d)$ of the conjugation by $s$ in
$E(d)$ and $C_{\sigma}$ denotes the restriction to $\Gamma$ of the conjugation by $\sigma$ in $Aut(G)$.

\medskip \noindent
With $(p,\pi)$ given, it follows from (\ref{symPair}) that $s$ is uniquely determined by $\sigma$. Indeed, assuming
the origin in $R^d$ to be the image by $p$ of a particular vertex $v_0\in G$, we have:

\begin{equation}\label{isometry}
 s(x)=Sx+t, 
\end{equation}

\begin{equation}\label{isometryBis}
 \mbox{with} \ S=C_s \ \mbox{when} \ R^d\equiv {\cal T}(R^d)\ \  \mbox{and} \ \  t=s(0)=p(\sigma v_0) 
\end{equation}

\medskip \noindent
Note that $s$ is a translation if and only if $C_s$ is the identity, that is, if and only if
$\sigma$ belongs to the centralizer ${\cal C}(\Gamma)$ of $\Gamma$ in $N(\Gamma)$:

$$ {\cal C}(\Gamma)=\{ \delta\in N(\Gamma)\ | \ \delta \gamma=\gamma \delta 
\ \mbox{for all} \ \gamma\in \Gamma \} $$

\medskip \noindent
Thus, the set of all symmetries of $(G,\Gamma,p,\pi)$ becomes a
group under composition. This is the {\em symmetry group} of the framework and will be denoted by
$\Sigma=\Sigma(G,\Gamma,p,\pi)\subset N(\Gamma)\subset Aut(G)$. The periodicity group $\Gamma$ is a normal
subgroup of $\Sigma$ and the injective homomorphism $\sigma \mapsto s$ is an extension of $\pi: \Gamma \rightarrow {\cal
T}(R^d)$ to $\Sigma \rightarrow E(d)$. For simplicity, this extension is also denoted by $\pi$.
Since $\pi(\Sigma)$ is a crystallographic group, we may refer to the framework $(G,\Gamma,p,\pi)$ as having {\em crystallographic symmetry} $\Sigma$.

\section{Symmetry constraints}

For a given $\sigma \in N(\Gamma)$, we may identify all periodic placements $(p,\pi)$ for
which $\sigma$ is a symmetry of the framework $(G,\Gamma,p,\pi)$. It will be convenient
to give this description in terms of parameters based on the following {\em choices}: a complete set of representatives $v_0,v_1,..., v_{n-1}$ for the vertex orbits of $\Gamma$ on $V$ and an
isomorphism $Z^d\approx \Gamma$. Later on,  a complete set of edge representatives of the form $e_{ij}=(v_i,v_j+\gamma_{ij})$ for the orbits of $\Gamma$ on $E$ will be implicated in obtaining equations for the length preservation of edges under deformation.

\medskip \noindent
Note that we allow the additive notation $\gamma v=v+\gamma$ for the action of
$\gamma \in  \Gamma$ on a vertex of the graph, which will facilitate writing the 
corresponding translation in a placement $(p,\pi)$ as $\pi(\gamma)p(v)=p(v)+\lambda$,
when the translation $\pi(\gamma)\in {\cal T}(R^d)$ has the formula $\pi(\gamma)(x)=x+\lambda$ and is identified with the translation vector $\lambda\in R^d$.

\medskip \noindent
Besides the frequent identification ${\cal T}(R^d)\equiv R^d$, which induces an inner product
on the group of translations, other routine {\bf conventions and notations} will be the following.
With the chosen isomorphism $Z^d\approx \Gamma$, the automorphisms of the free Abelian
group $\Gamma$ are given by $d\times d$ matrices with integer entries and determinant 
$\pm 1$ , that is $Aut(\Gamma)$ is identified with the {\em unimodular group} $GL(d,Z)$.
Moreover, by turning the free module $Z^d$ into the $d$-dimensional vector space

$$ R^d=Z^d\otimes_Z R\approx \Gamma \otimes_Z R $$

\noindent
unimodular transformations will be conceived as linear transformations. In particular,
conjugation by $\sigma \in N(\Gamma)$ gives an automorphism $C_{\sigma}\in Aut(\Gamma)$
and thereby a unimodular transformation which we denote by the same symbol
$C_{\sigma}\in GL(d,Z)$. Thus, we have a representation $N(\Gamma) \rightarrow GL(d,Z)$
with kernel ${\cal C}(\Gamma)$. 

\medskip \noindent
When given a periodic placement $(p,\pi)$ of $(G,\Gamma)$, $\pi$ gives an isomorphism

\begin{equation}\label{isom}
R^d\approx \Gamma \otimes_Z R \rightarrow {\cal T}(R^d)\equiv R^d
\end{equation}

\medskip \noindent
Thus, two bases are at play: the fixed Cartesian standard basis of $R^d\equiv {\cal T}(R^d)$
and the {\em lattice basis} which depends on $\pi: Z^d\approx \Gamma
\rightarrow {\cal T}(R^d)$. Sometimes, the coordinates based on the Cartesian basis are called
{\em geometric} and those based on the periodicity lattice are called {\em arithmetic}. 

\medskip
We have seen above that a necessary condition for $\sigma\in N(\Gamma)$ to be a symmetry
of the framework $(G,\Gamma,p,\pi)$ is that $C_{\sigma}\in GL(d,Z)$ should be orthogonal, when expressed in Cartesian coordinates. The remaining conditions refer to the translation
part $t$ in (\ref{isometry}), namely

\begin{equation}\label{t}
t=p(\sigma(v_i)) -S(p(v_i))
\end{equation}

\noindent
which must be the same for all $i=0,1,...,n-1$. It will be useful to express this in
arithmetic coordinates.  Let us assume that $\sigma$ acts on the vertex representatives
$v_0,v_1,...,v_{n-1}$ according to the formulae:

\begin{equation}\label{action}
\sigma(v_i)=v_{\sigma(i)}+\gamma_i
\end{equation}

\noindent
where $\sigma(i)$ is the index corresponding to the permutation effect of $\sigma$ on $V/\Gamma$ and $\gamma_i\in \Gamma$. Recall that the arithmetic or lattice coordines
are introduced through $\pi: \Gamma \rightarrow {\cal T}(R^d)\equiv R^d$ and the chosen
identification $Z^d\approx \Gamma$. Note that the periods  $\lambda_i\in \Gamma$ in (\ref{action}) will correspond with $n_i\in Z^d$ and these vectors with integer entries depend
only on $\sigma$ and the lattice identification $Z^d\approx \Gamma$. We allow here the same symbol $C_{\sigma}$ for the conjugation given by $\sigma$ on $\Gamma$, its expression as an
automorphism of $Z^d$ and its extension to $Z^d\otimes_Z R$.

\medskip  \noindent
With $e_k,\ 1\leq k\leq d$ the standard basis in
$Z^d$, we let $\Lambda_{\pi}=\Lambda\in GL(d,R)$ denote the  matrix with columns given by the lattice
basis $\pi(e_k)$. We define vector parameters $t_i$ by

\begin{equation}\label{representatives}
p(v_i)=\Lambda t_i
\end{equation}

\noindent
and obtain from (\ref{t}) and (\ref{action}) the following conditions:

\begin{equation}\label{t-part}
t_{\sigma(i)}+n_i-C_{\sigma}t_i=t_{\sigma(j)}+n_j-C_{\sigma}t_j  ,\ \ \mbox{for}\ \ 
0\leq i,j\leq n-1
\end{equation}

\medskip \noindent
Let us put $\omega_{\pi}=\omega=\Lambda^t \Lambda$ for the Gram matrix of the period lattice basis. Then, the orthogonality condition for  $C_{\sigma}\in GL(d,Z)$ becomes:

\begin{equation}\label{S-part}
C_{\sigma}^t\omega C_{\sigma}=\omega
\end{equation}

\medskip \noindent
In summary, we have

\begin{prop}\label{SymConstraints}
A graph automorphism $\sigma\in N(\Gamma)= Aut(G,\Gamma)$ is a symmetry of the $d$-periodic framework $(G,\Gamma,p,\pi)$ if and only if conditions (\ref{t-part}) and (\ref{S-part}) are satisfied.
\end{prop}

\medskip \noindent
Recall that the placement information $(p,\pi)$ enters in these equations through the 
parameters $t_i, \ 0\leq i\leq n-1$ and $\omega$ as described above. In subsequent sections, we shall elaborate on their role in describing symmetric periodic placements and symmetry preserving deformations.

\section{Parametrizations}

The deformation theory  developed in \cite{BS2} for periodic bar-and-joint frameworks in
$R^d$ emphasized the analogy with the traditional theory of finite linkages. In particular,
equivalent realizations resulting from isometries applied to any given framework were not 
immediately factored out. However, enumerative purposes or other concerns require the
quotient operation. In the finite case \cite{BS1}, Cayley-Menger matrices or equivalently, Gram matrices serve the purpose. In the periodic case, crystallography and  lattice theory have proven long ago the importance of the identification of the quotient  $O(d,R)\backslash GL(d,R)$
with the space of {\em positive definite quadratic forms} in $d$ variables, itself represented 
by the open cone $\Omega(d)$ of symmetric $d\times d$ matrices with positive eigenvalues \cite{CS, G, S1, W}.

\medskip \noindent
The parametrization used in the previous section follows this classical perspective. 
All the information about the lattice of periods $\pi(\Gamma)$, up to orthogonal
transformations, is contained in the symmetric matrix $\omega=\Lambda^t \Lambda\in 
\Omega(d)$, while the `shift vectors' $t_i$ indicate (relative to the lattice basis) the
placement of the vertex representatives $v_i$. By requesting that $t_0=0$, equivalence under translation is eliminated as well. This yields

\begin{prop}\label{param}
Let $(G,\Gamma)$ be a $d$-periodic graph. Then all periodic placements in $R^d$,
up to equivalence under the group of Euclidean isometries $E(d)$, are parametrized
by $(R^d)^{n-1}\times \Omega(d)$, which is an open set of $R^{dn+{d\choose 2}}$.
\end{prop}
  
\medskip \noindent
{\bf Remarks.}\  The vertex image sets of periodic placements are multilattices and this type of
configurations has been considered in different contexts. 
While not implicating an edge structure, the study of multilattices envisaged in \cite{P, PZ} 
is related to the kinematics of phase transitions in crystalline materials. 
When approached from the point of view of periodic sphere packings, as in \cite{S1, S2}, multilattices do acquire an edge structure from contacts between spheres. 
The resulting  {\em packing frameworks} are a very particular class of periodic frameworks. 
A study of homogeneous sphere packings in dimension three 
has been undertaken by W. Fischer, E. Koch and H. Sowa in a series of papers e.g. \cite{KSF}. 
The planar homogeneous case goes back to Niggli \cite{N1, N2}. See also \cite{F}.

\medskip
The bar-and-joint understanding of a framework brings in the (squared) length function
for edges and the notion of deformations \cite{BS2}. For a given $d$-periodic framework 
$(G,\Gamma,p,\pi)$ vertices become joints and edges become straight rigid bars between them.
It is enough to register the (squared) length of a complete set of representatives for $E/\Gamma$. As mentioned earlier, with vertex representatives $v_0,...,v_{n-1}$ already
chosen, we may select  $m$ edge representatives of the form $e_{ij}=(v_i,v_j+\gamma_{ij})$.
For expressing the squared length in both geometric and arithmetic coordinates, we let
$n_{ij}\in Z^d$ stand for the vector with $\Lambda n_{ij}$ equal to the translation vector of $\pi(\gamma_{ij})$. Then

\begin{equation}\label{SqLength}
 \ell(e_{ij})^2= || p(v_j+\gamma_{ij})-p(v_i) ||^2=(t_j+n_{ij}-t_i)^t\omega (t_j+n_{ij}-t_i)
\end{equation}

\medskip \noindent
This gives a {\em polynomial map} with cubic $(i\neq j)$ or linear $(i=j)$  components.

\begin{equation}\label{SqFunction}
\begin{array}{c} f: (R^d)^{n-1}\times \Omega(d) \rightarrow R^m \\
                           \  \\
 f(t_1,...,t_{n-1},\omega)= ((t_j+n_{ij}-t_i)^t\omega (t_j+n_{ij}-t_i))_{ij} \in R^m 
\end{array}
\end{equation}

\noindent
where the $m$ pairs of indices $ij$ correspond to the chosen representatives for 
the edge orbits $E/\Gamma$.

\medskip \noindent
The non-empty fibers of this map are {\em configuration spaces of frameworks} and the
connected component of a framework in its configuration space gives the {\em deformation space of the framework}.

\medskip 
For a simple comparison of this treatment with our presentation in \cite{BS2}, we offer
the following diagram.

\begin{equation}\label{OldNew}
 \begin{array}{clcl}
 (R^d)^n\times GL(d,R) & \stackrel{q}  \longrightarrow & E(d)\backslash (R^d)^n\times GL(d,R)  \\
                           \downarrow  & \   & \parallel   \\
                           R^m & \stackrel{f} \longleftarrow   & (R^d)^{n-1}\times \Omega(d) 
     \end{array}
\end{equation}

\noindent
$ (R^d)^n\times GL(d,R)$ is the parametrization used in \cite{BS2} for periodic placements
which include all isometric replicas of all frameworks. Hence, the fibers of the left vertical arrow
are {\em realization spaces} for weighted periodic graphs $(G,\Gamma,\ell)$, that is, 
periodic graphs with prescribed lengths for their edges. When isometries are factored out,
we obtain, as stated above in Proposition~\ref{param}, the parameter space $(R^d)^{n-1}\times \Omega(d)$ with the bottom map (\ref{SqFunction}). With explicit formulae, we have the following description.

\medskip  
$ (R^d)^n\times GL(d,R)$ parametrizes periodic placements $(p,\pi)$ by recording the positions
of the $n$ vertex representatives and the basis of the lattice of periods $\pi(\Gamma)$, that is,
$(p(v_0),...,p(v_{n-1}),\Lambda)$. The left action of the isometry group $E(d)$ on these parameters is given by

\begin{equation}\label{leftAction}
u(p(v_0),...,p(v_{n-1}),\Lambda)=(u\circ p(v_0),...,u\circ p(v_{n-1}),U\Lambda)
\end{equation}

\noindent
for an isometry $u(x)=Ux+t$, with $U\in O(d,R)$ and $t\in R^d$. The {\em quotient map} $q$
works by the formula

\begin{equation}\label{quotientMap}
\begin{array}{c} q(p(v_0),...,p(v_{n-1}),\Lambda)=(t_1,...,t_{n-1},\omega)= \\
                           \  \\
=(\Lambda^{-1}(p(v_1)-p(v_0)),...,\Lambda^{-1}(p(v_{n-1})-p(v_0)),\Lambda^t\Lambda)  \end{array}
\end{equation}

\noindent
The left vertical arrow is the composition $f\circ q$.

\medskip 
A direct enumerative consequence of the current presentation will be an upper bound for the
number of distinct possible configurations of a {\em minimally rigid periodic graph} with generic edge lenght prescriptions. Recall from \cite{BS2, BS3} that minimally rigid periodic graphs have
$m=dn+{d\choose 2}$ edge orbits. In the generic case, the corresponding edge length
constraints are independent. There can be no more than ${d+1\choose 2}$ linear constraints $(i=j)$ among them, because all linear constraints affect only $\omega$. Since the polynomial map (\ref{SqFunction}) can be extented to
{\em complex projective} coordinates in $P_m(C)$, we infer from B\'{e}zout's theorem the following bound.

\begin{prop}\label{bound}
Let  $(G,\Gamma)$ be a minimally rigid $d$-periodic graph with $n=|V/\Gamma|$ and $m=|E/\Gamma|=nd+{d\choose 2}$. Let $\mu$ be the number of cubic edge constraints
$(i\neq j)$ in (\ref{SqFunction}). Then $dn-d\leq \mu \leq m$ and $(G,\Gamma)$  has at most $3^{\mu}$ non-congruent configurations in $R^d$ for a generic prescription of edge lengths.
\end{prop}

\medskip \noindent
{\bf Remark.} This upper bound result is analogous to the one  obtained in \cite{BS1}
for finite minimally rigid graphs. 

\section{Actions and representations}

We may now return to symmetry considerations and elaborate on the {\em affine} nature
of the symmetry constraints obtained in Proposition~\ref{SymConstraints} in terms of 
placement parameters $(t_1,...,t_{n-1},\omega)$. 

\medskip
Let us recall that, by definition, $\sigma\in Aut(G,\Gamma)$ becomes a symmetry of a placement $(p,\pi)$ when the effect of $\sigma$ on the periodic graph $(G,\Gamma)$ is {\em reproduced} by the effect of an isometry $s\in E(d)$
on the image of the graph determined by $p(V)$. In other words, the placements $(p,\pi)$ and
$(p\circ \sigma, \pi\circ C_{\sigma})$ must be equivalent under the action of $E(d)$. Hence,
in parameters $(t_1,...,t_{n-1},\omega)$, we must have one and the same point. 
By Proposition~\ref{SymConstraints}, the {\em fixed point locus} of such an action by $\sigma$
is given by an {\em affine linear subvariety} and this fact leads to the obvious expectation that 
the action itself is expressed by an affine map in the parameters $(t_1,...,t_{n-1},\omega)$. 

\medskip \noindent
This is indeed the case, as ensuing computations will confirm. In order to discuss the 
action of $Aut(G,\Gamma)$ on placements and the quotient parameter space 
$(R^d)^{n-1}\times \Omega(d)\subset R^{dn+{d\choose 2}}$ as a {\em left action}, we adopt the following convention.

\begin{definition}\label{LeftAction}
Let $\sigma\in Aut(G,\Gamma)$ be an automorphism of a $d$-periodic graph $(G,\Gamma)$.
The left action of $Aut(G,\Gamma)$  on periodic placements in $R^d$ is defined by
the formula:

\begin{equation}\label{onLeft}
\sigma (p,\pi)=(p\circ \sigma^{-1}, \pi\circ C_{\sigma^{-1}})
\end{equation}
 
\end{definition}

\begin{theorem}\label{AffineRep}
When expressed in parameters $(t_1,...,t_{n-1},\omega)$, the action (\ref{onLeft}) corresponds
with an affine representation

\begin{equation}\label{Rep}
A: Aut(G,\Gamma)\rightarrow Aff(dn+{d\choose 2})
\end{equation}

\noindent
which factors through $Aut(G,\Gamma)/\Gamma=Aut(G/\Gamma)$.

\medskip \noindent
For any subgroup $\Gamma \subset \Sigma \subset Aut(G,\Gamma)$, the periodic placements 
of $(G,\Gamma)$ with crystallographyc symmetry $\Sigma$ are parametrized by the fixed locus 
of $A(\Sigma)$ in $(R^d)^{n-1}\times \Omega(d)$, that is

\begin{equation}\label{locus}
{\cal F}(\Sigma)=\{ x\in (R^d)^{n-1}\times \Omega(d) \ :\ 
A(\sigma)x=x, \ \mbox{for all}\ \sigma\in \Sigma  \}
\end{equation}

\noindent
The locus with full symmetry  ${\cal F}(Aut(G,\Gamma))$ is not empty.
\end{theorem}

\medskip \noindent
{\em Proof:}\ \  For the computation in coordinates $(t_1,...,t_{n-1},\omega)$, let us put

$$ \sigma(t_1,...,t_{n-1},\omega)=(\tilde{t}_1,...,\tilde{t}_{n-1},\tilde{\omega}) $$

\noindent
and recall that $t_0=\tilde{t}_0=0$. We have $\tilde{\omega}=\tilde{\Lambda}^t\tilde{\Lambda}$, with 
$\tilde{\Lambda}=\Lambda C_{\sigma^{-1}}$, hence 

\begin{equation}\label{explicit}
\tilde{\omega}=(C_{\sigma^{-1}})^t\omega C_{\sigma^{-1}}
\end{equation}

\medskip \noindent
Recall also that ${\sigma^{-1}}$ induces a permutation on $\{ 0,...,n-1\}$ by its effect on $V/\Gamma$. By (\ref{action}) and (\ref{quotientMap}), we find 
$\sigma^{-1}(v_j)=v_{\sigma^{-1}(j)}-C_{\sigma^{-1}}\gamma_{\sigma^{-1}(j)}$ and then

\begin{equation}\label{explicitBis}
\tilde{t}_j=C_{\sigma}(t_{\sigma^{-1}(j)}-t_{\sigma^{-1}(0)})+
(n_{\sigma^{-1}(0)}-n_{\sigma^{-1}(j)})
\end{equation}

\medskip \noindent
Formulae (\ref{explicit}) and (\ref{explicitBis}) give the explicit form of the action of $\sigma$,
which is linear in the components of $\omega$ and affine in the components of $t_j$, $j=1,...,n-1$.

\medskip \noindent
The resulting homomorphism (\ref{Rep}) is obviously trivial on $\Gamma$. Since $Aut(G,\Gamma)/\Gamma=Aut(G/\Gamma)$ is finite, the image group must have at least one fixed point (the barycenter of an orbit). 

\begin{cor}\label{posets}
There is an inclusion reversing correspondence $\Sigma \mapsto {\cal F}(\Sigma)$ between
subgroups $\Gamma\subset \Sigma\subset Aut(G,\Gamma)$ and a finite system of 
non-empty affine linear sections of $(R^d)^{n-1}\times \Omega(d)$ which parametrize periodic 
placements with a specified crystallographic symmetry.
\end{cor}

\medskip \noindent
It may be observed that this approach obtains periodic placements for $(G,\Gamma)$ with
full symmetry $Aut(G,\Gamma)$ realized by corresponding crystallographic groups, {\em without
recourse to a minimizing principle}. Other methods for proving the
existence of placements with higher symmetry rely explicitly on some `energy functional' minimization technique for finding `barycentric placements' \cite{DF} or a harmonic `standard placement' \cite{KS}.

\section{Relaxing or refining symmetry}

Up to this point, our considerations have focused on a given $d$-periodic graph $(G,\Gamma)$
with framework placements $(G,\Gamma,p,\pi)$ in $R^d$. However, various problems may
require a relaxation $\tilde{\Gamma}\subset \Gamma$ or a refinement $\Gamma\subset \hat{\Gamma}$ of the periodicity group. 
With the perspective gained in the preceding sections, we may introduce
the following definition.

\begin{definition}\label{commensurate}
Let $G=(V,E)$ be an infinite graph and let $\Gamma_1, \Gamma_2\subset Aut(G)$ be free
Abelian groups of rank $d$ such that the corresponding $d$-periodic graphs $(G,\Gamma_i)$
admit periodic presentations in $R^d$. Then, $\Gamma_1$ and $\Gamma_2$ are called commensurate when $\Gamma_1\cap \Gamma_2$ is of finite index in both groups $\Gamma_1$ and $\Gamma_2$.
\end{definition}

\medskip \noindent
Let us observe the effect of {\em relaxing} periodicity from $\Gamma$ to a subgroup $\tilde{\Gamma}\subset \Gamma$  of index $k$. We select a complete set of representatives
$\nu_0=0,\nu_1,...,\nu_{k-1}$ for $\Gamma/\tilde{\Gamma}$. Then, representatives for $V/\tilde{\Gamma}$ are given by $v_i+\nu_j$, with $0\leq i\leq n-1$ and $0\leq j\leq k-1$.

\medskip \noindent
When each periodicity group is
identified with $Z^d$, the inclusion of $\tilde{\Gamma}$ in $\Gamma$ corresponds to an
invertible matrix with integer entries $M$ with $det(M)=k$. Of course, the $\Gamma$-periodic placements of $(G,\Gamma)$ are contained in the $\tilde{\Gamma}$-periodic placements of
$(G,\tilde{\Gamma})$, and for a placement $(p,\pi)$  this inclusion takes the form:

\begin{equation}\label{relax}
((t_i)_i,\omega) \mapsto ((\tilde{t}_{ij})_{ij},\tilde{\omega}), \ \ \mbox{with}\ t_0=0=\tilde{t}_{00}
\end{equation}

\noindent
With $\pi(\nu_j)=m_j$ as translation vectors, we have $p(v_i+\nu_j)=p(v_i)+m_j$.  By (\ref{quotientMap}) and its counterpart for $\tilde{\Gamma}$, we find:

\begin{equation}\label{inclusion}
\tilde{\Lambda}=\Lambda M, \ \ \mbox{hence}\ 
\tilde{\omega}=\tilde{\Lambda}^t\Lambda=M^t\omega M
\end{equation}

\begin{equation}\label{inclusionBis}
\tilde{t}_{ij}=M^{-1}t_i+M^{-1}\Lambda^{-1}m_j
\end{equation}

\medskip \noindent
We already know that, from the perspective of $(G,\tilde{\Gamma})$, the placements with
`higher' symmetry $\Gamma$ are parametrized by an affine linear section and the above computation confirms the expected fact that we have an affine inclusion map which 
identifies the parameter space for periodic placements of $(G,\Gamma)$ with this affine
linear section. It follows that {\em the affine and convexity structure} of the
parameter space $(R^d)^{n-1}\times \Omega(d)\subset R^{dn+{d\choose 2}}$ for periodic
placements of $(G,\Gamma)$ is preserved when relaxing the lattice.

\medskip  
Recall that a {\em crystallographic group} in dimension $d$ is a {\em discrete} subgroup of isometries $K\subset E(d)$, with a {\em compact} quotient $E(d)/K$. Bieberbach showed that
the subgroup of translations in $K$, that is $L=K\cap {\cal T}(R^d)$, must be a lattice of rank
$d$ and is uniquely determined as the maximal free Abelian normal subgroup of $K$. Moreover, $K/L$ is a finite group. Thus, if $\Sigma\subset Aut(G)$ is isomorphic with a crystallographic group $K\subset E(d)$, we may refer to the free Abelian normal subgroup $\Gamma\subset \Sigma$ corresponding to $L\subset K$, and form the $d$-periodic graph $(G,\Gamma)$. 

\medskip \noindent
We have seen above that, when $(G,\Gamma)$ allows periodic placements in $R^d$, some of
them, namely those parametrized by ${\cal F}(\Sigma)$, will have Euclidean symmetries 
given by some crystallographic group isomorphic with $\Sigma$ and $K$ (which must be, 
according to another Bieberbach theorem, an affine conjugate of $K$). 

\medskip \noindent
Under these circumstances, we may refer to $\Sigma\subset Aut(G)$ as a {\em crystallographic subgroup} of $Aut(G)$ and use the pair notation $(G,\Sigma)$ for the graph $G$ {\em with the specified crystallographic symmetry} $\Sigma$. It is also understood that, up to Euclidean isometry, the placements of $(G,\Sigma)$ are those parametrized by ${\cal F}(\Sigma)$. We note
that ${\cal F}(\Sigma)$ can be determined in the placement parameter space of any periodic
graph $(G,\tilde{\Gamma})$ with $\tilde{\Gamma}\subset \Gamma$ of finite index and stable
under conjugation by $\Sigma$. This determination amouts to solving a linear system of
equations with integer coefficients of the form (\ref{t-part}) and (\ref{S-part}), corresponding to a finite set of transformations $\sigma\in \Sigma$ which provide generators for $\Sigma/\tilde{\Gamma}$. Verifications entirely similar to those performed above show that
{\em the affine and convexity structure} of ${\cal F}(\Sigma)$ is the same for all choices of
$\tilde{\Gamma}$.

\medskip 
The commensurability equivalence relation extends as follows.

\begin{definition}\label{commensurateBis}
Two crystallographic subgroups $\Sigma_1, \Sigma_2\subset Aut(G)$ are called commensurate
when $\Sigma_1\cap \Sigma_2$ is of finite index in both $\Sigma_1$ and $\Sigma_2$.
\end{definition}

\medskip \noindent
{\bf Remarks.}\ It would be enough to assume $\Sigma_1\cap \Sigma_2$ of finite index in one of the groups, but we prefer the symmetric formulation. Clearly, in this case, the intersection
$\Sigma_1\cap \Sigma_2$ is itself a crystalographic subgroup of $Aut(G)$. Subgroups of finite index in crystallographic groups can be found by simple procedures \cite{Se}.

\medskip  
{\em Relaxing or refining the symmetry} of a framework are associated with certain variations  within the framework's commensurability equivalence class.
This language seems favorable for addressing geometric aspects of {\em displacive phase
transitions} in crystalline materials \cite{Bu, D}. The vertices of the infinite graph $G$ may
serve as labels for a subfamily or all of the atoms in some idealized crystal, with edges marking
bonds. Under variations of temperature or pressure, the same material may have phases with 
different crystallographic symmetry. Displacive phase transitions involve no bond rupture, hence the graph remains the same. When two phases have commensurate symmetry groups $\Sigma_1, \Sigma_2$, our approach gives the simplest geometrical common ground for a passage, namely
${\cal F}(\Sigma_1\cap \Sigma_2)$, which contains both ${\cal F}(\Sigma_1)$ and 
${\cal F}(\Sigma_2)$ as affine linear sections.

\medskip \noindent
Of course, as a guiding scenario, this has been formulated long ago. The new insight, at least at
the geometric level, is that the symmetry preserving loci have a simple affine structure and description. However, non-linearity resurfaces when bonds are assumed to maintain
their length. This is considered in the next section.

\section{Symmetry preserving deformations}

When returning to the edge squared length function (\ref{SqFunction}), a first simple remark is that $f$ is $Aut(G,\Gamma)$ equivariant. To see this, it is convenient to write $R^m$ as the 
space $R^{E/\Gamma}$ of real-valued functions on $E/\Gamma$. Then the left action of
$Aut(G,\Gamma)$ is simply defined as 

\begin{equation}\label{Eaction}
\sigma(\phi)=\phi\circ \sigma^{-1}, \ \ \mbox{for} \ \phi\in R^{E/\Gamma}
\end{equation}

\noindent
with the action on $E$ and on $E/\Gamma$  denoted by the same symbol.
Then, one  easily verifies that

\begin{equation}\label{equivariant}
f(\sigma(p,\pi))=\sigma(f(p,\pi))
\end{equation}

\medskip \noindent
When given a framework $(G,\Gamma,p,\pi)$ with crystallographic symmetry $\Sigma\subset
N(\Gamma)=Aut(G,\Gamma)$, and we want to consider {\em only} deformations which 
preserve this symmetry (and all edge lengths), we have to restrict $f$ to ${\cal F}(\Sigma)$
and consider the fiber of $f(p,\pi)$. Since $\Sigma$ acts trivially on ${\cal F}(\Sigma)$, the
image by $f$ must consist of points invariant under $\Sigma$, that is, $f$ factors through
$R^{E/\Sigma}$.

\medskip \noindent
{\bf Remark.}\ Strictly speaking, we should write $\Sigma\backslash E$ for the quotient by
an action on the left. Expecting no harm, we continue with $E/\Sigma$ and note that
$E/\Gamma\rightarrow E/\Sigma$ induces $R^{E/\Sigma}\rightarrow R^{E/\Gamma}$.

\medskip \noindent
It follows that the edge length control for frameworks with crystallographic symmetry $\Sigma$
is given by a map which we allow to be denoted by the same symbol

\begin{equation}\label{SqFunctionBis}
f: {\cal F}(\Sigma) \rightarrow R^{E/\Sigma}
\end{equation}

\medskip \noindent
Thus, after appropriate rank computations, we obtain a setting entirely analogous to the 
basic case $\Sigma=\Gamma$.

\vspace{0.2in}
 
   Ciprian S. Borcea\\
   Department of Mathematics\\
   Rider University\\
   Lawrenceville, NJ 08648, USA

\vspace{0.2in} 

   Ileana Streinu\\
   Department of Computer Science\\
   Smith College\\ 
   Northampton, MA 01063, USA

\end{document}